\documentclass[reqno,11pt]{amsart}
\usepackage{latexsym,amssymb,amsfonts,amsmath,amsthm}
\usepackage[usenames]{color}
\usepackage{eucal}
\usepackage{mathrsfs}
\usepackage{setspace}

\newtheorem{theorem}{Theorem}[section]
\newtheorem{corollary}[theorem]{Corollary}
\theoremstyle{remark}
\newtheorem{remark}[theorem]{Remark}
\numberwithin{equation}{section}
\newcommand{\Real}{\mathbb R}
\newcommand{\Z}{\mathbb Z}
\newcommand{\eps}{\varepsilon}

\newcommand{\F}{\mathcal{F}}

\newcommand{\C}{\mathcal{C}}
\newcommand{\one}[1]{\mathbf{1}_{\{#1\}}}
\newcommand{\oneset}[1]{\mathbf{1}_{#1}}
\newcommand{\PP}{\mathbb{P}}

\newcommand{\E}{\mathbb{E}}

\newcommand{\iunit}{\mathbf{i}}

\newcommand{\Law}{\mathscr{L}}

\renewcommand{\theenumi}{A\arabic{enumi}}

\begin{document}

\title[]{Compound Poisson approximation for triangular arrays with  application to threshold estimation}%
\author{P. Chigansky}%
\address{Department of Statistics,
The Hebrew University,
Mount Scopus, Jerusalem 91905,
Israel}
\email{pchiga@mscc.huji.ac.il}

\author{F.C. Klebaner}
\address{ School of Mathematical Sciences,
Monash University Vic 3800,
Australia}
\email{fima.klebaner@monash.edu}

%\author{et al}
\thanks{ Research supported by the  Australian Research Council Grant
DP0988483. P.Chigansky was supported by ISF grant 314/09.}


\date{21, October, 2011}%
%\dedicatory{}%
%\commby{}%
% ----------------------------------------------------------------
\begin{abstract}
We prove  weak convergence of triangular arrays to the
compound Poisson limit using Tihomirov's method.  The result is applied to statistical estimation of the threshold  parameter
in autoregressive models.
\end{abstract}
\maketitle

\section{Introduction and main result}

This paper is concerned with  weak convergence of sums over triangular arrays with certain dependence structure
to the compound Poisson distribution. It is motivated by the threshold estimation problem, described in details in Section \ref{a-sec}.
We consider triangular arrays of random variables  $Y_{n,j}$, $j=1,...,n$, $n\in \mathbb{N}$  with rows, adapted to a filtration
$(\F_j)$, $j\in \mathbb{N}$. $Y_{n,j}$'s are asymptotically negligible and satisfy a weak dependence (mixing)  condition made precise  by the following assumptions.

\medskip

\begin{enumerate}

\item\label{a1} there is a constant $C_1>0$, such that
\[
\PP(Y_{n,j}\ne 0)\le \frac {C_1}{n}, \quad \text{and}\quad \E|Y_{n,j}|\le \frac {C_1}{n}, \quad j=1,...,n
\]
and
\[
\quad \E|Y_{n,j}|\one{Y_{n,i}\ne 0}\le \left(\frac {C_1}{n}\right)^2, \quad i\ne j.
\]

\item\label{a2} there is an integer $\ell\ge 1$, such that
\[
\Big|\E \big(Y_{n,j}\big|\F_{i}\big)-\E Y_{n,j}\Big| \le \frac{C_1}{n} \alpha(j-i),\quad i\le j-\ell
\]
where  $\alpha(n)\ge 0$ is a decreasing sequence with $\lim_{n\to\infty}\alpha(n)=0$

\medskip

\item\label{a3} for  a measurable function $|v(x)|\le 1$, $x\in \Real^{n-j+1}$
\[
\Big|\E \big(v(Y_{n,j},...,Y_{n,n})\big|\F_{i}\big)-\E v(Y_{n,j},...,Y_{n,n})\Big| \le \alpha(j-i),\quad i<j
\]
\end{enumerate}

\medskip

The following condition on the individual characteristic functions $\phi_{n,j}(t)=\E e^{\iunit t Y_{n,j}}$
together with the above assumptions, will assure convergence of the sums
\[
S_n = \sum_{j=1}^n Y_{n,j}, \quad n\in \mathbb{N},
\]
to the compound Poisson law
(hereafter we shall abbreviate $\dot \varphi(t)=\frac d {dt} \varphi(t)$, etc.):

\medskip
\begin{enumerate}
\setcounter{enumi}{3}
\item\label{a4}  There exists  a characteristic function $\varphi(t)$ and positive constants  $C_2$ and $\mu$  such that
\[
\Big|\dot \phi_{n,j}(t) -n^{-1} \mu \dot \varphi(t)\Big|\le C_2 n^{-2}, \quad t\in \Real.
\]
\end{enumerate}
\medskip

Note that the mixing in \ref{a2} and \ref{a3} can be arbitrarily weak. Further assumptions on  the rate of convergence of  $\alpha(k)$ to zero, such
as:
\medskip

\begin{enumerate}
\setcounter{enumi}{4}
\item\label{a5}  $\alpha(k)\le C_3 r^k$ for some $r\in (0,1)$ and $C_3>0$.
\end{enumerate}

\medskip
\noindent
allow to obtain rates of convergence in an appropriate metric.
Below we shall work with the L\'evy distance, defined
for a pair of distribution functions $F$ and $G$ by (see e.g. \cite{GS02})
\[
L(G,F)=\inf\big\{h>0:G(x-h)-h\le F(x)\le G(x+h)+h, \
\forall \ x\big\}.
\]

Our main result is the following:

\begin{theorem}\label{thm}
Let $Y_{n,j}$, $j=1,...,n$, $n\in \mathbb{N}$ be a triangular array of
random variables, whose rows are adapted to a filtration $(\F_j)$, $j\in \mathbb{N}$ and satisfy the assumptions \ref{a1}-\ref{a4}.
 Then
\begin{equation}
\label{wc}
S_n = \sum_{j=1}^n Y_{n,j}\xrightarrow[n\to\infty]{d}S,
\end{equation}
where $S$ has the compound Poisson distribution, with intensity $\mu$ and i.i.d. jumps with characteristic
function $\varphi(t)$.

Moreover, if the assumption \ref{a5} holds then
 there is a constant $C>0$, such that for all $n$ large
enough,
\begin{equation}
\label{Lbnd}
L\big(\Law(S_n),\Law(S)\big) \le C  n^{-1/2}\log n,
\end{equation}
where $L\big(\Law(S_n),\Law(S)\big)$ is the L\'{e}vy distance between the distribution functions of $S_n$ and $S$.

\end{theorem}

\begin{remark}
Both the constant $C$ and the smallest $n$ for which \eqref{Lbnd} holds, can be found explicitly in terms of
the $C_i$'s and $\alpha(\cdot)$, mentioned in the assumptions above. Also bounds on the L\'{e}vy distance can be obtained similarly for e.g.
polynomially decreasing $\alpha(\cdot)$, by replacing $b\log n$ with $n^\delta$ for some $\delta>0$ in the proof of
Theorem 1.1 and optimizing the right hand side of the corresponding inequality, analogous to \eqref{Zol} below.
\end{remark}

In  application to threshold estimation, $Y_{n,j}$ is  derived from an autoregressive stationary process $X_j$, generated by the recursion
\begin{equation}
\label{X}
X_j = h(X_{j-1})+\eps_j, \quad j\ge 1,
\end{equation}
where $h(\cdot)$ is a given measurable function and $(\eps_j)$ is a sequence of i.i.d. random variables, with continuous positive probability
density $q(\cdot)$. As explained in Section \ref{a-sec}, in this context
\[
Y_{n,j}:= f(\eps_j) \one{X_{j-1}\in B_n},
\]
where     $B_n:= [0, 1/n]$, $f(\cdot)$ is a measurable function and
\[
S_n := \sum_{j=1}^n f(\eps_j) \one{X_{j-1}\in B_n}.
\]
Theorem \ref{thm} implies that under appropriate conditions, $S_n$ converges weakly to the compound  Poisson random variable with i.i.d. jumps, distributed as $f(\eps_1)$, and the intensity
$\mu:=p(0)$, where $p(\cdot)$ is the unique invariant density of $(X_j)$.

Somewhat surprisingly, we were not able to find in the literature a general result, from which this limit could be deduced.
In this regard, one naturally thinks of Stein's method or martingale convergence results.
Stein's method appears to be particularly well suited to the compound Poisson distribution with {\em integer} valued jumps
(see e.g. \cite{BC01}, \cite{AGG}). The results such as \cite{B80}, \cite{S86, S86c}, \cite{BV10} or \cite{R03} come close, but apparently do not quite fit our setting.

In the particular case, when $\E f(\eps_j) = 0$, $S_n$ becomes a sum over the array of  martingale differences
$Y_{n,j}:= f(\eps_j) \one{X_{j-1}\in B_n}$, $j=1,...,n$ with the quadratic variation sequence
\[
V_{n,m} = \sum_{j=1}^m \one{X_{j-1}\in B_n}\E f^2(\eps_1), \quad m=1,...,n.
\]
A typical martingale limit result such as e.g. \cite{BE71} or Theorem 2.27 Ch. VIII \S 2c in \cite{JSh2ed}
requires that $V_{n,n}$ converges in probability.
However  in our case  $V_{n,n}$   converges only in distribution  (to a Poisson random variable), but not in probability
(since e.g. $V_{n,n}$ is uniformly integrable, but is not a Cauchy sequence in $L_1$).  It is known that
$S_n$ may have a different limit or no limit at all, if the convergence in probability of quadratic variation
is replaced with convergence in distribution (see \cite{ACF81} and the references therein), so that the martingale
results also do not appear applicable\footnotemark.

\footnotetext{in this connection, it is interesting to note, that in the analogous continuous time setting, the quadratic variation
does converge in probability, essentially due to the continuity of the sample paths, see \cite{K11}}

The objective of this paper is to give a  proof of Theorem \ref{thm}, using Tihomirov's method from \cite{T80}. Originally
applied to CLT in the dependent case, it turns to be remarkably suitable to the setting under consideration. Before proceeding to the proof in
Section \ref{p-sec}, we shall discuss in more details the application, in which the aforementioned convergence arises.

\section{Application to threshold estimation}\label{a-sec}
  Suppose one observes a sample $X^n=(X_1,...,X_n)$
from a threshold autoregressive (TAR) time series, generated by the recursion
\begin{equation}
\label{TAR}
X_j = g_+(X_{j-1}) \one{X_{j-1}\ge \theta}+g_-(X_{j-1})\one{X_{j-1}<\theta}+ \eps_j, \quad j\in \Z_+,
\end{equation}
where $g_+(\cdot)$ and $g_-(\cdot)$ are known functions and $(\eps_j)$ is a sequence of i.i.d. random variables with
known probability density $q(\cdot)$.  The unknown {\em threshold} parameter $\theta$, taking values in an open interval $\Theta:=(a,b)\subset \Real$, is
to be estimated from the sample $X^n$. TAR models, such as \eqref{TAR}, have been the subject of extensive research in statistics and econometrics
(see e.g. \cite{Tong83} and the recent surveys \cite{Tong11}, \cite{H11}, \cite{ChK10}).

From the statistical analysis point of view, this estimation problem classifies as ``singular'', since
the corresponding likelihood function
\begin{equation}
\label{L}
L_n(X^n;\theta)= \prod_{j=1}^n q\Big(X_j - g_+(X_{j-1}) \one{X_{j-1}\ge \theta}-g_-(X_{j-1})\one{X_{j-1}<\theta}\Big)
\end{equation}
is discontinuous in $\theta$. Typically in such problems, the sequence of the Bayes estimators
\[
\tilde \theta_n := \frac{\int_\Theta \theta L_n(X^n;\theta)\pi(\theta)d\theta}{\int_\Theta  L_n(X^n;\theta)\pi(\theta)d\theta}, \quad n\ge 1
\]
is  asymptotically efficient in the  minimax sense for an arbitrary continuous prior density  $\pi(\cdot)$ (see \cite{IH}).
The asymptotic distribution of these estimators is determined by the weak limit of the likelihood ratios as follows.
Let $\theta_0\in \Theta$ be the true unknown value of the parameter and $r_n$ an increasing sequence of numbers.  The change of variables
$u=r_n(\theta-\theta_0)\in r_n(\Theta-\theta_0)=:\mathbb{U}_n$ gives
\[
r_n(\tilde \theta_n -\theta) =
\frac{\int_{\mathbb{U}_n} u Z_n(u)\pi(\theta_0+u/r_n)du}{\int_{\mathbb{U}_n} Z_n(u)\pi(\theta_0+u/r_n)du},
\]
where $Z_n(u)$, $n\ge 1$ are the rescaled likelihood ratios
\[
Z_n(u) = \frac{L_n(X^n;\theta_0+u/r_n)}{L_n(X^n;\theta_0)}, \quad u\in \mathbb{U}_n.
\]
If $r_n$ can be chosen so that $Z_n(u)$, $u\in \Real$ converges weakly to a random process $Z(u)$, $u\in \Real$ in an appropriate topology, then
\begin{equation}
\label{asym}
r_n(\tilde \theta_n -\theta) \xrightarrow[n\to\infty]{d} \frac{\int_{\Real} u Z(u)du}{\int_{\Real} Z(u)du},
\end{equation}
holds (a comprehensive account of this approach can be found in \cite{IH}).

For the likelihoods as in \eqref{L}, a simple calculation (see eq. (4) in \cite{ChK}) reveals that
\begin{equation}
\label{lnZnu}
\log Z_n(u) =
\sum_{j=1}^n \one{X_{j-1}\in B_n}\log \frac{q\big(\eps_j+\delta(X_{j-1})\big)}{q\big(\eps_j\big)}
,\quad u\ge 0
\end{equation}
where $B_n := [\theta_0,\theta_0+u/n]$ and $\delta(x):=g_+(x)-g_-(x)$, and a similar expression is obtained for $u<0$.
It can be shown that \eqref{asym} indeed holds with  $r_n:=n$, if $(X_j)$ is a sufficiently fast
mixing with the unique invariant probability density $p(x;\theta_0)$, and
the sequence $\log Z_n(u)$ converges weakly to the compound Poisson process
\begin{equation}
\label{lnZu}\log Z(u):=
\begin{cases}
\sum_{j=1}^{\Pi^+(u)}
\log \frac{q\big(\eps^+_j+\delta(\theta_0)\big)}{q\big(\eps^+_j\big)} & u\ge 0\\
\sum_{j=1}^{\Pi^-(-u)}
\log \frac{q\big(\eps^-_j-\delta(\theta_0)\big)}{q\big(\eps^-_j\big)} & u< 0
\end{cases},
\end{equation}
where $(\eps^{\pm}_j)$ are i.i.d. copies of $\eps_1$ and $\Pi^{+}(u)$ and $\Pi^{+}(u)$ are independent Poisson processes with the same intensity
$p(\theta_0;\theta_0)$.

The rate $r_n=n$ and the Poisson behavior  is typical for discontinuous likelihoods (see e.g. Ch. 5, \cite{IH}).
For the linear TAR  model, i.e. when $g_{\pm}(x)=\rho_{\pm}x$ with constants $\rho_-\ne \rho_+$, this asymptotic appeared  in \cite{Ch93} and the aforementioned generalization is taken from \cite{ChK}.

One particularly interesting ingredient in the proof, which is the main focus of this article, is the convergence of the finite
dimensional distributions of $Z_n(u)$ to those of $Z(u)$. In its  prototypical form, the problem can be restated as follows.
Consider the stationary Markov sequence $(X_j)$, generated by the recursion \eqref{X} and let (cf. \eqref{lnZnu})
\begin{equation}
\label{Sn}
S_n := \sum_{j=1}^n f(\eps_j) \one{X_{j-1}\in B_n},
\end{equation}
where $B_n:= [0, 1/n]$ and $f(\cdot)$ is a measurable function. It is required to show that,
the sums $(S_n)$ converge weakly to the compound  Poisson random variable with i.i.d. jumps, distributed as $f(\eps_1)$, and the intensity
$p(0)$, where $p(\cdot)$ is the unique invariant density of $(X_j)$.

This convergence is not hard to prove using the blocks technique: $S_n$ is partitioned into, say, $n^{1/2}$  blocks of $n^{1/2}$
consecutive summands, $n^{1/4}$ of which are discarded. Removing total of $n^{1/2}\cdot n^{1/4}$ out of $n$ terms in the sum does not alter its limit,
but the residual blocks become nearly independent, if the mixing is fast enough. Moreover, a single event $\{X_{j}\in B_n\}$
occurs within each block with probability of order $n^{-1/2}$ and hence the sum over approximately independent $n^{1/2}$ blocks
yields the claimed compound Poisson behavior. This approach dates back to at least \cite{M73} in the Poisson case, and the details for the
compound Poisson setting can be found in \cite{ChK}.

\medskip

An alternative proof now can be given by applying Theorem \ref{thm}:

\begin{corollary}\label{cor}
Let $(X_j)$ be  defined by \eqref{X} and $S_n$ by \eqref{Sn}.
Assume that

\medskip

\begin{enumerate}
\renewcommand{\theenumi}{\roman{enumi}}
\item\label{i} $\eps_1$ has positive Lipschitz continuous  bounded probability density $q(x)$, $x\in \Real$ with the finite first absolute moment
$\int_\Real |x|q(x)dx<\infty$

\medskip

\item\label{ii}  for some $r\in (0,1)$ and $C>0$,
\[
|h(x)|\le r|x|, \quad \forall\; |x|\ge C
\]

\medskip

\item\label{iii}  $\E |f(\eps_1)|<\infty$ and for some constant $C'$,
\[
\sup_{z,x\in [0,n^{-1}]} |f\big(z-h(x)\big)|\le C'
\]
for all $n$ large enough.
\end{enumerate}
\medskip

Then the Markov process $(X_j)$ has   unique invariant density $p(x)$, $x\in \Real$, which is positive, Lipschitz continuous  and bounded;
for stationary  $(X_j)$, the sums $(S_n)$ converge weakly to the compound Poisson random variable with intensity $p(0)$ and i.i.d. jumps
with the same distribution as $f(\eps_1)$.
\end{corollary}

\begin{remark}
The Corollary \ref{cor} verifies the weak convergence of the one-dimensional distributions of the processes $\log Z_n(u)$ from \eqref{lnZnu}
to those of $\log Z(u)$, $u\in \Real$ defined in \eqref{lnZu}. The convergence of finite dimensional distributions of higher orders can be
treated along the same lines. The limit \eqref{asym} then follows from the tightness of the sequence of processes $\log Z_n(u)$
(see \cite{ChK} for further details).
\end{remark}
\begin{remark}
The assumption \ref{iii} holds if e.g. $f(\cdot)$ and $h(\cdot)$ are continuous at $0$.
\end{remark}

\begin{proof}
Under the assumptions \ref{i} and \ref{ii}, the standard ergodic theory of Markov chains (see e.g. Theorem 16.0.2 in \cite{MT09}) implies that $(X_j)$
is irreducible, aperiodic  and positive recurrent Markov chain with the unique invariant measure. Due to the additive structure of the recursion \eqref{X}, the invariant measure
has density $p(\cdot)$, which is positive and  continuous with the same Lipschitz constant $L_q$ as the density $q(\cdot)$ and
$\|p\|_\infty\le \|q\|_\infty:=\sup_{x\in \Real}q(x)$.
Moreover, $(X_j)$ is   geometrically mixing, i.e.  there exist positive constants $R$ and $\rho<1$, such that for any
measurable function $|g(x)|\le 1$
\begin{equation}
\label{mix}
\left|\E\big(g(X_j)|\F_i\big) - \int_\Real g(x)p(x)dx\right|\le R \rho^{j-i}, \quad j>i,
\end{equation}\
where $\F_j=\sigma\{\eps_i, i\le j\}$.
Define $Y_{n,j}:= f(\eps_j)\one{X_{j-1}\in B_n}$, then
\[
\E|Y_{n,j}|= \E|f(\eps_j)|\PP(X_{j-1}\in B_n) =
 \E\big|f(\eps_1)\big| \int_0^{1/n} p(x)dx\le \E|f(\eps_1)| \|q\|_\infty n^{-1},
\]
and similarly
\[
\PP(Y_{n,j}\ne 0) \le  \PP(X_{j-1}\in B_n)\le \|q\|_\infty n^{-1}.
\]
Further, for  $i<j-1$,
\begin{align*}
\E\big|Y_{n,i}  \one{Y_{n,j}\ne 0}\big|=\; &
\E\big|f(\eps_i)|\one{X_{i-1}\in B_n}  \PP(Y_{n,j}\ne 0|\F_i)\le \\
&\E\big|f(\eps_i)|\one{X_{i-1}\in B_n}  \PP(X_{j-1}\in B_n|\F_i)\le  \\
&\E\big|f(\eps_i)|\one{X_{i-1}\in B_n} \|q\|_\infty n^{-1}\le
\E\big|f(\eps_1)| \|q\|^2_\infty n^{-2},
\end{align*}
and
\begin{multline*}
\E\big|  \one{Y_{n,i}\ne 0}Y_{n,j}\big|=
\E  \one{Y_{n,i}\ne 0}\E\big(|Y_{n,j}|\big|\F_i\big)\le \\
\E  \one{X_{i-1}\in B_n}\E\big(\big|f(\eps_j)\big|\one{X_{j-1}\in B_n}\big|\F_i\big)\le
\E \big|f(\eps_1)\big|\|q\|^2_\infty n^{-2}.
\end{multline*}
Similarly,
\begin{align*}
\E|Y_{n,j-1}| & \one{Y_{n,j}\ne 0} \le
\E \one{X_{j-2}\in B_n}  |f(\eps_{j-1})| \one{X_{j-1}\in B_n}= \\
&
\E \one{X_{j-2}\in B_n} \E\Big( |f(\eps_{j-1})| \one{h(X_{j-2})+\eps_{j-1}\in B_n}|\F_{j-2}\Big)=\\
&
\E \one{X_{j-2}\in B_n} \int_\Real |f(y)| \one{h(X_{j-2})+y\in B_n}q(y)dy=\\
&
\E \one{X_{j-2}\in B_n} \int_0^{n^{-1}} \big|f\big(z-h(X_{j-2})\big)\big| q\big(z-h(X_{j-2})\big)dz\le\\
&
\|q\|_\infty\E \one{X_{j-2}\in B_n} n^{-1}\sup_{z,x\in [0,n^{-1}]} |f\big(z-h(x)\big)|\le
\|q\|^2_\infty n^{-2}C'
\end{align*}
and
\[
\E  \one{Y_{n,j-1}\ne 0} |Y_{n,j}|\le
\E  \one{X_{j-2}\in B_n} |f(\eps_j)|\one{X_{j-1}\in B_n}\le
\E  \big|f(\eps_1)\big|  \|q\|^2_\infty n^{-2}.
\]
Hence \ref{a1} is satisfied for all $n$ large enough with
\[
C_1: =  \Big(\|q\|^2_\infty\vee 1\Big) \Big(\E\big|f(\eps_1)\big| \vee C'\vee 1\Big).
\]
Further, by the Markov property,
\begin{multline*}
\E (Y_{n,j}|\F_{j-2}) =
\E\Big(\one{X_{j-1}\in B_n}\E \big(f(\eps_j)|\F_{j-1}\big)\big|\F_{j-2}\Big)=\\
\E f(\eps_1)\PP\big(X_{j-1}\in B_n\big|X_{j-2}\big)=: H(X_{j-2}),
\end{multline*}
and
\[
\big|H(X_{j-2})\big|\le \big|\E f(\eps_1)\big|\|q\|_\infty n^{-1}\le C_1 n^{-1}.
\]
Hence by \eqref{mix}, for $i<j-1$,
\begin{multline*}
\Big|\E \big(Y_{n,j}\big|\F_{i}\big)-\E Y_{n,j}\Big| =
\Big|\E\Big(\E \big(Y_{n,j}\big|\F_{j-2}\big)|\F_i\Big)-\E\E \big(Y_{n,j}\big|\F_{j-2}\big)\Big|=\\
\Big|\E\big(H(X_{j-2})|\F_i\big)-\E H(X_{j-2})\Big| \le
C_1 n^{-1}R \rho^{j-i-2},
\end{multline*}
and \ref{a2} holds with $\ell=2$ and
\begin{equation}\label{mixAR1}
\alpha(k):=R \rho^{k-2}.
\end{equation}
 The assumption \ref{a3} is checked similarly.
Finally,
\begin{multline*}
\dot \phi_{n,j}(t) = \E \iunit f(\eps_j)\one{X_{j-1}\in B_n} e^{\iunit t f(\eps_j)\one{X_{j-1}\in B_n}} =
\E \iunit f(\eps_j)\one{X_{j-1}\in B_n} e^{\iunit t f(\eps_j)}= \\
\E \one{X_{j-1}\in B_n}\E\Big(\iunit f(\eps_j) e^{\iunit t f(\eps_j)}\big|\F_{j-1}\Big) =
\PP(X_1\in B_n) \dot \varphi(t),
\end{multline*}
where $\varphi(t)=\E e^{\iunit t f(\eps_1)}$ and interchanging derivative and the expectation is valid by the
dominated convergence and \ref{iii}.

Since the invariant density is Lipschitz, it follows that
\[
\Big|\dot\phi_{n,j}(t)-p(0)\frac 1 n \dot \varphi(t)\Big|\le L_q n^{-2},
\]
which verifies \ref{a4} and the claim now follows from Theorem \ref{thm}.
In fact, the assumption \ref{a5} holds by virtue of \eqref{mixAR1} and  the L\'evy distance to the limit distribution
converges at the rate, claimed in \eqref{Lbnd}.
\end{proof}

\section{Proof of Theorem \ref{thm}}\label{p-sec}

Tihomirov's approach \cite{T80} is applicable, when the characteristic function of the limit distribution uniquely solves  an ordinary
differential equation.  Roughly, the idea is then to show that the characteristic functions of the prelimit distributions satisfy the same
equation in the limit.

The characteristic function of the compound Poisson distribution with intensity $\mu$ and characteristic function of the jumps $\varphi(t)$
is given by
\[
\psi(t) = e^{\mu(\varphi(t)-1)},\quad t\in \Real
\]
which solves uniquely the initial value problem
\[
\dot \psi(t) =  \mu \dot \varphi(t) \psi(t), \quad \psi(0)=1, \quad t\in \Real.
\]

Since $\E|S_n|<\infty$, the characteristic function $\psi_n(t):= \E e^{\iunit t S_n}$ is
continuously differentiable and   $\Delta_n(t):= \psi(t)-\psi_n(t)$ satisfies
\[
\dot \Delta_n(t) = \mu \dot \varphi(t)\Delta_n(t) + r_n(t), \quad t\in \Real,
\]
subject to $\Delta_n(0)=0$, where $r_n(t):= \mu \dot\varphi(t)\psi_n(t)-\dot \psi_n(t)$. Solving for $\Delta_n(t)$ gives
\begin{equation}
\label{exps}
\Delta_n(t) = \int_0^t \exp\Big(\mu \big(\varphi(t)-\varphi(s)\big)\Big) r_n(s)ds, \quad t\ge 0.
\end{equation}
As we show below, for any constant $b>0$, such that $b\log n$ is a positive integer,
\begin{equation}
\label{showme}
|r_n(t)|\le
C_2 n^{-1} +
3C_1 \alpha(b\log n)+
8C^2_1 b \frac {\log n}{n}, \quad t\ge 0
\end{equation}
and, since $|\varphi(t)|\le 1$, it follows from \eqref{exps} that
\begin{equation}\label{D}
|\Delta_n(t)|\le e^{2\mu}\int_0^t  |r_n(s)|ds\le
 e^{2\mu}\left(
C_2 n^{-1} +
3C_1 \alpha(b\log n)+
8C^2_1 b \frac{\log n}{n}
\right)t, \quad t\ge 0.
\end{equation}
Similar bound holds for $t<0$ and the claimed weak limit \eqref{wc} follows, once we check \eqref{showme}.
To this end, we have
\begin{align*}
& \dot \psi_n(t) :=  \frac d {dt}\E e^{\iunit tS_n} =  \E \frac d {dt} \exp \Big(\iunit t \sum_{j=1}^n Y_{n,j}\Big)=
\sum_{k=1}^n \E \iunit   Y_{n,k} \exp \Big(\iunit t \sum_{j=1}^n Y_{n,j}\Big)=\\
&
 \sum_{k=1}^n \E\iunit   Y_{n,k} e^{\iunit t Y_{n,k}}\exp \Big(\iunit t \sum_{j\ne k} Y_{n,j}\Big)=
 \sum_{k=1}^n \E\iunit   Y_{n,k} e^{\iunit t Y_{n,k}}\exp \Big(\iunit t \sum_{|j- k|>b\log n} Y_{n,j}\Big)
+ \\
& \sum_{k=1}^n \E\iunit   Y_{n,k} e^{\iunit t Y_{n,k}}\bigg(
\exp \Big(\iunit t \sum_{j\ne k} Y_{n,j}\Big)
-
\exp \Big(\iunit t \sum_{|j- k|>b\log n} Y_{n,j}\Big)
\bigg) := J_1+J_2
\end{align*}
where we used $\E|S_n|<\infty$ and the dominated convergence  to interchange the derivative and the expectation.
Note that $\big|e^{\iunit x}-e^{\iunit (x+y)}\big|\le 2\one{y\ne 0}$ for any $x,y\in \Real$, and hence
by  the assumption \ref{a1}
\begin{align}\label{J2}
& |J_2| \le
\sum_{k=1}^n \E  |Y_{n,k}|\bigg|
\exp \Big(\iunit t \sum_{j\ne k} Y_{n,j}\Big)
-
\exp \Big(\iunit t \sum_{|j- k|>b\log n} Y_{n,j}\Big)
\bigg|\le \\
\nonumber
& 2\sum_{k=1}^n \E  |Y_{n,k}|\oneset{ \left\{\sum_{|j- k|\le b\log n, j\ne k} Y_{n,j}\ne 0\right\}}\le
2\sum_{k=1}^n \E  |Y_{n,k}| \sum_{|j- k|\le b\log n, j\ne k}\one{ Y_{n,j}\ne 0}\le
4C_1^2 b\frac {\log n}{n}.
\end{align}
Further,  by the triangle inequality
\begin{align*}
&
\bigg|\E\iunit   Y_{n,k} e^{\iunit t Y_{n,k}}\exp \bigg(\iunit t \sum_{|j- k|>b\log n} Y_{n,j}\bigg) - \dot \phi_{n,k}(t)\psi_n(t)\bigg|=\\
& \bigg|\E\iunit   Y_{n,k} e^{\iunit t Y_{n,k}}\exp \bigg(\iunit t \sum_{|j- k|>b \log n} Y_{n,j}\bigg) - \E\iunit   Y_{n,k} e^{\iunit t Y_{n,k}}
\E \exp \bigg(\iunit t \sum_{j=1}^n Y_{n,j}\bigg)\bigg|\le \\
&
\bigg|\E   Y_{n,k} e^{\iunit t Y_{n,k}}\exp \bigg(\iunit t \sum_{|j- k|> b \log n} Y_{n,j}\bigg) - \E  Y_{n,k} e^{\iunit t Y_{n,k}}
\E \exp \bigg(\iunit t \sum_{|j- k|>b \log n} Y_{n,j}\bigg) \bigg|
+\\
&
\Big|\E   Y_{n,k} e^{\iunit t Y_{n,k}}\Big|\bigg|\E\exp \bigg(\iunit t \sum_{|j- k|>b\log n} Y_{n,j}\bigg) -
\E \exp \bigg(\iunit t \sum_{j=1}^n Y_{n,j}\bigg) \bigg|=: J_3+J_4.
\end{align*}
Similarly to \eqref{J2}, we have
\begin{multline*}
|J_4|\le
\E   |Y_{n,k}|\E \bigg|\exp \bigg(\iunit t \sum_{|j- k|>b \log n} Y_{n,j}\bigg) -
\exp \bigg(\iunit t \sum_{j=1}^n Y_{n,j}\bigg) \bigg|\le \\
 2\E   |Y_{n,k}|\E\oneset{\left\{
\sum_{|j- k|\le b\log n} Y_{n,j}\ne 0
\right\}}\le
2\E   |Y_{n,k}|
\sum_{|j- k|\le b\log n} \PP(Y_{n,j}\ne 0)\le 4 C^2_1 b\frac{\log n}{n^2}.
\end{multline*}
For brevity, define
\[
U := \exp \bigg(\iunit t \sum_{j<k-b\log n} Y_{n,j}\bigg), \quad
V :=    Y_{n,k} e^{\iunit t Y_{n,k}}, \quad
W :=\exp \bigg(\iunit t \sum_{j>k+b\log n} Y_{n,j}\bigg).
\]
By the triangle inequality,
\begin{multline*}
  \big|
\E UVW-\E V \E UW
\big|\le
\big|
\E UVW-\E UV\E W
\big|+ \\
\big|
\E UV \E W-
\E U\E V \E W
\big|
+
\big|
\E U\E V \E W
-\E V \E U W
\big|.
\end{multline*}
Since $U$ and $V$ are $\F_k$-measurable, $|U|\le 1$, $|W|\le 1$ and $\E|V|\le C_1n^{-1}$, \ref{a3} implies
\[
\big|
\E UVW-\E UV\E W
\big|\le \E|UV|\big|
\E(W|\F_k)-\E W
\big|\le C_1 n^{-1}  \alpha(b\log n),
\]
and, since $U$ is measurable with respect to $\F_{k-b\log n}$,
\begin{multline*}
\big|
\E U\E V \E W
-\E V \E U W
\big|\le  |\E V|\E|U|
\big|
\E W
- \E ( W|\F_{k-b\log n})
\big|\le C_1 n^{-1} \alpha(2b\log n).
\end{multline*}
Further, by \ref{a2} for $b\log n\ge \ell$,
\[
\big|
\E UV \E W-
\E U\E V \E W
\big|\le
|\E W|\E|U|
\big|
\E(V|\F_{k-b\log n}) -
\E V
\big|\le C_1 n^{-1} \alpha(b\log n).
\]
Hence
\[
|J_3|=
\Big|
\E UVW-\E V \E UW
\Big|\le 3C_1 n^{-1} \alpha(b\log n),
\]
and consequently, by \ref{a4}
\begin{multline*}
\Big|J_1-\mu \dot \varphi(t)\psi_n(t)\Big|=
\bigg|
\sum_{k=1}^n \E\iunit   Y_{n,k} e^{\iunit t Y_{n,k}}\exp \Big(\iunit t \sum_{|j- k|>b\log n} Y_{n,j}\Big)-\mu \dot \varphi(t)\psi_n(t)\bigg|\le
  \\C_2 n^{-1} +
3C_1 \alpha(b\log n)+
4C^2_1 b\frac {\log n}{n}.
\end{multline*}
Assembling all parts together, we obtain \eqref{showme}.

The  bound \eqref{Lbnd} for the L\'{e}vy metric is obtained by means of Zolotorev's inequality   \cite{Zol},
\[
L\Big(\Law(S_n),\Law(S)\Big)\le \frac 1 \pi \int_0^T \frac{|\psi_n(t)-\psi(t)|}{t}dt +2e \frac{\log T}{T}, \quad T>1.3,
\]
which in view of the bound in \eqref{D} gives
\begin{equation}
\label{Zol}
L\Big(\Law(S_n),\Law(S)\Big)\le
\frac  {e^{2\mu}}\pi
\Big(
C_2 n^{-1} +
3C_1 \alpha(b\log n)+
2C^2_1 b\frac{\log n}{n}
\Big)T
 +2e \frac{\log T}{T}.
\end{equation}
If $\alpha(k)$ decays geometrically as in \ref{a5}, the bound  \eqref{Lbnd} is obtained by choosing $T=n^{1/2}$ and $b\ge \frac 1{\log 1/r}$.
\qed

\begin{remark}
The rate in \eqref{Lbnd} is not as sharp as the one, obtained by Tihomirov in \cite{T80} in the CLT case.
Apparently, the deficiency originates in the specific form of the compound Poisson characteristic function $\psi(t) = e^{\mu(\varphi(t)-1)}$,
which does not vanish as $t\to\infty$. More specifically, the integration kernel $K(s,t):=e^{\mu (\varphi(t)-\varphi(s))}$ in
\eqref{exps} does not decay when $t$ is fixed and $s$ decreases, which contributes the linear growth in $t$ of the right hand side of
\eqref{D} and the corresponding linear growth in $T$ in \eqref{Zol}. In the Gaussian case, this kernel has the form $K(s,t):=e^{s^2/4-t^2/4}$
(see eq. (3.25) page 809 in \cite{T80}), which yields  better balance between growth in $t$ and the decrease in $n$.
It seems that in the compound Poisson setting under consideration the rate cannot be essentially improved within the
framework of Tihomirov's method.
\end{remark}

\subsection*{Acknowledgement} The authors are grateful to Y. Kutoyants and R. Liptser for
the enlightening discussions on the subject of this article.
We also appreciate referee's suggestion, which led to the ultimate
improvement of the rate in \eqref{Lbnd}.
% ----------------------------------------------------------------

%\bibliographystyle{amsplain}
%\bibliography{compound}

\def\cprime{$'$}
\providecommand{\bysame}{\leavevmode\hbox to3em{\hrulefill}\thinspace}
\providecommand{\MR}{\relax\ifhmode\unskip\space\fi MR }
\providecommand{\MRhref}[2]{%
  \href{http://www.ams.org/mathscinet-getitem?mr=#1}{#2}
}
\providecommand{\href}[2]{#2}

\end{document}